\documentclass[11pt]{amsart}
\usepackage[T1]{fontenc}
\usepackage{amssymb}
\usepackage{amsmath}
\usepackage{fancyhdr}
\usepackage[british]{babel}
\usepackage{geometry}
\usepackage{enumitem}
\usepackage{algpseudocode}
\usepackage{dsfont}
\usepackage{centernot}
\usepackage{xstring}
\usepackage{colortbl}
\usepackage[section]{algorithm}
\usepackage{graphicx}
\usepackage[font={footnotesize}]{caption}
\usepackage[usenames,dvipsnames,table]{xcolor}
\usepackage{pstricks,tikz}
\usepackage[h]{esvect}
\usepackage[
		bookmarksopen=true,
		bookmarksopenlevel=1,
		colorlinks=true,
		linkcolor=darkblue,
        linktoc=page,
		citecolor=darkblue,
]{hyperref}

\title[]{Euler tours in hypergraphs}
\date{\today}
\author{Stefan Glock}
\author{Felix Joos}
\author{Daniela K\"uhn}
\author{Deryk Osthus}

\address{School of Mathematics, University of Birmingham,
Edgbaston, Birmingham, B15 2TT, United Kingdom}
\email{ [s.glock, f.joos, d.kuhn, d.osthus]@bham.ac.uk}

\thanks{The research leading to these results was partially supported by the European Research Council under the European Union's Seventh Framework Programme (FP/2007--2013) / ERC Grant 306349 (S.~Glock and D.~Osthus), and by the Deutsche Forschungsgemeinschaft (DFG, German Research Foundation) -- 339933727~(F.~Joos). The research was also partially supported by the EPSRC, grant no. EP/N019504/1, and by the Royal Society and the Wolfson Foundation (D.~K\"uhn).}

\newtheorem{theorem}{Theorem}

\newtheorem{lemma}[theorem]{Lemma}

\theoremstyle{definition}

\newtheorem{conj}[theorem]{Conjecture}

\newtheoremstyle{claimstyle}{5pt}{5pt}{\em}{5pt}{\em}{:}{5pt}{}
\theoremstyle{claimstyle}

\newtheoremstyle{stepstyle}{10pt}{5pt}{\em}{0pt}{\em}{:}{5pt}{}
\theoremstyle{stepstyle}

\definecolor{darkblue}{rgb}{0,0,0.5}

\def\noproof{{\unskip\nobreak\hfill\penalty50\hskip2em\hbox{}\nobreak\hfill%
       $\square$\parfillskip=0pt\finalhyphendemerits=0\par}\goodbreak}
\def\endproof{\noproof\bigskip}

\newdimen\margin
\def\textno#1&#2\par{
   \margin=\hsize
   \advance\margin by -4\parindent
          \setbox1=\hbox{\sl#1}
   \ifdim\wd1 < \margin
      $$\box1\eqno#2$$
   \else
      \bigbreak
      \hbox to \hsize{\indent$\vcenter{\advance\hsize by -3\parindent
      \it\noindent#1}\hfil#2$}
      \bigbreak
   \fi}

\def\proof{\removelastskip\penalty55\medskip\noindent\setcounter{claim}{0}\setcounter{step}{0}{\bf Proof. }} % in each main proof, claim counter set back

\def\lateproof#1{\removelastskip\penalty55\medskip\noindent\setcounter{claim}{0}\setcounter{step}{0}{\bf Proof of #1. }} % in each main proof, claim counter set back

\begin{document}

\newcommand{\new}[1]{\textcolor{red}{#1}}
\def\COMMENT#1{}
\def\TASK#1{}

\newcommand{\todo}[1]{\begin{center}\textbf{to do:} #1 \end{center}}

\def\eps{{\varepsilon}}
\newcommand{\ex}{\mathbb{E}}
\newcommand{\pr}{\mathbb{P}}
\newcommand{\cB}{\mathcal{B}}
\newcommand{\cA}{\mathcal{A}}
\newcommand{\cE}{\mathcal{E}}
\newcommand{\cS}{\mathcal{S}}
\newcommand{\cF}{\mathcal{F}}
\newcommand{\cG}{\mathcal{G}}
\newcommand{\bL}{\mathbb{L}}
\newcommand{\bF}{\mathbb{F}}
\newcommand{\bZ}{\mathbb{Z}}
\newcommand{\cH}{\mathcal{H}}
\newcommand{\cC}{\mathcal{C}}
\newcommand{\cM}{\mathcal{M}}
\newcommand{\bN}{\mathbb{N}}
\newcommand{\bR}{\mathbb{R}}
\def\O{\mathcal{O}}
\newcommand{\cP}{\mathcal{P}}
\newcommand{\cQ}{\mathcal{Q}}
\newcommand{\cR}{\mathcal{R}}
\newcommand{\cJ}{\mathcal{J}}
\newcommand{\cL}{\mathcal{L}}
\newcommand{\cK}{\mathcal{K}}
\newcommand{\cD}{\mathcal{D}}
\newcommand{\cI}{\mathcal{I}}
\newcommand{\cV}{\mathcal{V}}
\newcommand{\cT}{\mathcal{T}}
\newcommand{\cU}{\mathcal{U}}
\newcommand{\cW}{\mathcal{W}}
\newcommand{\cX}{\mathcal{X}}
\newcommand{\cZ}{\mathcal{Z}}
\newcommand{\1}{{\bf 1}_{n\not\equiv \delta}}
\newcommand{\eul}{{\rm e}}
\newcommand{\Erd}{Erd\H{o}s}
\newcommand{\cupdot}{\mathbin{\mathaccent\cdot\cup}}
\newcommand{\whp}{whp }
\newcommand{\bX}{\mathcal{X}}
\newcommand{\bV}{\mathcal{V}}

\newcommand{\doublesquig}{%
  \mathrel{%
    \vcenter{\offinterlineskip
      \ialign{##\cr$\rightsquigarrow$\cr\noalign{\kern-1.5pt}$\rightsquigarrow$\cr}%
    }%
  }%
}

\newcommand{\defn}{\emph}

\newcommand\restrict[1]{\raisebox{-.5ex}{$|$}_{#1}}

\newcommand{\prob}[1]{\mathrm{\mathbb{P}}\left[#1\right]}
\newcommand{\expn}[1]{\mathrm{\mathbb{E}}\left[#1\right]}
\def\gnp{G_{n,p}}
\def\G{\mathcal{G}}
\def\lflr{\left\lfloor}
\def\rflr{\right\rfloor}
\def\lcl{\left\lceil}
\def\rcl{\right\rceil}

\newcommand{\qbinom}[2]{\binom{#1}{#2}_{\!q}}
\newcommand{\binomdim}[2]{\binom{#1}{#2}_{\!\dim}}

\newcommand{\grass}{\mathrm{Gr}}

\newcommand{\brackets}[1]{\left(#1\right)}
\def\sm{\setminus}
\newcommand{\Set}[1]{\{#1\}}
\newcommand{\set}[2]{\{#1\,:\;#2\}}
\newcommand{\krq}[2]{K^{(#1)}_{#2}}
\newcommand{\ind}[1]{$\mathbf{S}(#1)$}
\newcommand{\indcov}[1]{$(\#)_{#1}$}
\def\In{\subseteq}

\begin{abstract}  \noindent
We show that a quasirandom $k$-uniform hypergraph $G$ has a tight Euler tour subject to the necessary condition that $k$ divides all vertex degrees.
The case when $G$ is complete confirms a conjecture of Chung, Diaconis and Graham from 1989 on the existence of universal cycles for the $k$-subsets of an $n$-set.
\end{abstract}

\maketitle

\section{Introduction}

Finding an \defn{Euler tour} in a graph is a problem as old as graph theory itself: Euler's negative resolution of the Seven Bridges of K\"onigsberg problem in 1736 is widely considered the first theorem in graph theory. Euler observed that if a (multi-)graph contains a closed walk which traverses every edge exactly once, then all the vertex degrees are even. He also stated that every connected graph with only even vertex degrees contains such a walk, which was later proved by Hierholzer and Wiener.

There are several ways of generalising the concept of paths/cycles, and similarly Euler trails/tours, to hypergraphs. Not least due to its connection to universal cycles, we focus in this paper on the so-called `tight' regime. We will discuss related notions in Section~\ref{sec:rel results}.

\subsection{Universal cycles}
Let $[n]$ denote the set $\Set{1,\dots,n}$ and $\binom{[n]}{k}$ the set of all $k$-subsets of~$[n]$.  A \defn{universal cycle for $[n]$}
is a cyclic sequence with $\binom{n}{k}$ elements,
each of which is from $[n]$, such that every $k$ consecutive elements are distinct and every element of $\binom{[n]}{k}$ appears exactly once consecutively (but in an arbitrary order).
For example, $1234524135$ is a universal cycle for~$\binom{[5]}{2}$. The study of these objects was initiated by Chung, Diaconis and Graham~\cite{CDG:92} in a paper where they define universal cycles for various combinatorial structures (see Section~\ref{sec:rel results}). 

Observe that the number of $k$-subsets in $[n]$ that contain a particular element is $\binom{n-1}{k-1}$
and every element in a cyclic sequence (of length at least $k+1$) appears in exactly $k$ sets of $k$ consecutive elements.
Hence if a universal cycle for $\binom{[n]}{k}$ exists, then $k$ divides $\binom{n-1}{k-1}$, or equivalently, $n$ divides $\binom{n}{k}$. In 1989, Chung, Diaconis and Graham conjectured that the converse should also be true, at least if $n$ is sufficiently large, and offered $\$100$ for the resolution of this problem.

\begin{conj}[Chung, Diaconis, Graham~\cite{CDG:89,CDG:92}] \label{conj:CDG}
For every $k\in \bN$, there exists $n_0\in \bN$ such that for all $n\ge n_0$, there exists a universal cycle for $\binom{[n]}{k}$ whenever $k$ divides $\binom{n-1}{k-1}$.
\end{conj}

It is easy to see that Conjecture~\ref{conj:CDG} is true for $k=2$. 
Numerous partial results have been obtained. In particular, Jackson proved the conjecture for $k=3$~\cite{jackson:93} and for $k\in \Set{4,5}$ (unpublished), and Hurlbert~\cite{hurlbert:94} confirmed the cases $k\in\Set{3,4,6}$ if $n$ and $k$ are coprime (see also~\cite{LTT:81}).
Various approximate versions of Conjecture~\ref{conj:CDG} have been obtained in~\cite{blackburn:12,CHHM:09,DL:16,LTT:81}.
We prove Conjecture~\ref{conj:CDG} in a strong form by showing the existence of tight Euler tours in `typical' $k$-graphs.

\subsection{Tight Euler tours in typical hypergraphs}
Given a $k$-graph $G$ (i.e.~a $k$-uniform hypergraph~$G$), a sequence of vertices $\cW=x_1x_2\dots x_\ell$ is a \defn{(tight self-avoiding) walk in $G$} if $\Set{x_i,x_{i+1},\dots ,x_{i+k-1}}\in E(G)$ for all $i\in [\ell-k+1]$, and no edge of $G$ appears more than once in this way. Similarly, we say that $\cW$ is a \defn{closed walk} if $\Set{x_i,x_{i+1},\dots ,x_{i+k-1}}\in E(G)$ for all $i\in [\ell]$, with indices modulo~$\ell$, and no edge of $G$ appears more than once in this way. We let $E(\cW)$ denote the set of edges appearing in~$\cW$. An \defn{Euler tour of $G$} is a closed walk $\cW$ in~$G$ with $E(\cW)=E(G)$, and an \defn{Euler trail of $G$} is a walk $\cW$ in~$G$ with $E(\cW)=E(G)$.\COMMENT{could also require that the start and end $(k-1)$-sets are disjoint} Clearly, a universal cycle for $\binom{[n]}{k}$ is equivalent to an Euler tour of the complete $n$-vertex $k$-graph~$K^k_n$.

The problem of deciding whether a given $3$-graph has an Euler tour has been shown to be NP-complete~\cite{LNR:17}. Thus, when $k>2$, there is probably no simple characterisation of $k$-graphs having an Euler tour. However, we show that for `typical' $k$-graphs, the existence of an Euler tour hinges only on a simple divisibility condition.

A $k$-graph $G$ on $n$ vertices is called \defn{$(c,h,p)$-typical} if for any set $A$ of $(k-1)$-subsets of $V(G)$ with $|A|\le h$, we have $||\bigcap_{S\in A}N_G(S)|-p^{|A|}n|\leq cp^{|A|}n$, where $N_G(S)$ denotes the \defn{neighbourhood of $S$}, i.e.~the set of all vertices which together with $S$ form an edge.
Note that this is what one would expect in a random $n$-vertex $k$-graph in which every edge appears independently with probability~$p$. It is easy to see that the complete $k$-graph $K^k_n$ is $(hk/n,h,1)$-typical. Thus, the following more general result implies Conjecture~\ref{conj:CDG}.

\begin{theorem}\label{thm:typical}
For all $k\in \bN$ and $p\in(0,1]$,
there exist $c>0$ and $h,n_0\in\bN$ such that the following holds: Let $G$ be a $(c,h,p)$-typical $k$-graph on at least $n_0$ vertices with all vertex degrees divisible by~$k$.
Then $G$ has a tight Euler tour.
\end{theorem}
Clearly, the condition that all vertex degrees are divisible by $k$ is necessary for the existence of a tight Euler tour. 
Instead of an Euler tour, we can also easily obtain a tight Euler trail (see end of Section~\ref{sec:proof}).

We briefly sketch the strategy of our proof.
Let us first consider graphs.
We first find a closed walk in $G$ that contains all vertices. 
Afterwards we decompose the remainder of $G$ into small cycles and insert them into the closed walk to obtain an Euler tour.
For $k$-graphs, this can be done as follows.
In the first step, we find a `spanning' walk $\cW$ in~$G$, where spanning means that every ordered $(k-1)$-set of vertices appears at least once consecutively in the vertex sequence of~$\cW$. For this, we show that a self-avoiding random walk yields such a walk~$\cW$ (after an appropriate number of steps) with high probability. This step will use only a small fraction of the edges of~$G$. We then extend $\cW$ to a closed walk~$\cW'$. 
Subsequently, we remove $E(\cW')$ from $G$ and decompose the remaining $k$-graph into tight cycles using recent results of Glock, K\"uhn, Lo and Osthus~\cite{GKLO:17} (which imply the existence if $F$-designs). 
Each such cycle can be incorporated into~$\cW'$, which finally yields a tight Euler tour.

\subsection{Related research and open questions} \label{sec:rel results}
The most prominent example of universal cycles are de Bruijn cycles. A \defn{de Bruijn cycle of order~$k$} is a binary cyclic sequence in which every binary sequence of length~$k$ appears as a subsequence (of consecutive terms) exactly once. 
Chung, Diaconis and Graham~\cite{CDG:92} extended this notion to various other combinatorial objects, for instance permutations (see~\cite{johnson:09}) and partitions of an~$n$-set. 
The general idea is that a universal cycle for a set $S$ of combinatorial objects is a cyclic sequence which contains a `representation' of every element of~$S$ exactly once as a subsequence of consecutive terms. Due to their rich symmetry, such structures have found many applications, for instance in cryptography, computer graphics, database theory, digital fault testing and neural decoding. A common and natural approach to find universal cycles is via \defn{transition graphs}. Suppose that every element of $S$ has a unique representation as a sequence of length~$k$. One can then define a directed graph $G_S$ with vertex set $S$ where there is an arc from $(x_1,\dots,x_k)$ to $(y_1,\dots,y_k)$ if and only if $y_{i}=x_{i+1}$ for all $i\in[k-1]$. With this terminology, a universal cycle for $S$ corresponds to a directed Hamilton cycle in~$G_S$. One obstacle to finding universal cycles for $\binom{[n]}{k}$, which was noted in~\cite{CDG:92}, is that it is not even possible to define such a transition graph (since each $k$-set is represented by several sequences).

Rather than seeking an Euler tour in a $k$-graph~$G$, an alternative way to cover all edges of~$G$ is to ask for a Hamilton decomposition,
 i.e.~to ask for a collection of edge-disjoint Hamilton cycles in~$G$ such that every edge is contained in exactly one such Hamilton cycle. 
In 1892, Walecki showed that the complete (2-)graph~$K_n$ has a Hamilton decomposition whenever $n$ is odd.
As mentioned before, there are several natural definitions of paths/cycles in hypergraphs. One of the earliest such concepts was introduced by Berge. A \defn{Berge cycle} consists of a cyclic alternating sequence $v_1 e_1 v_2 e_2\dots v_\ell e_\ell$ of distinct vertices and edges such that $v_i,v_{i+1} \in e_i$ for all $i \in [\ell ]$.
(Here $v_{\ell+1} := v_1$ and the edges $e_i$ are also allowed to contain vertices outside $\Set{v_1,\dots,v_\ell}$.)
For $n \ge 100$, it is shown in \cite{KO:14b} that $K^k_n$ has a decomposition into Berge Hamilton cycles if and only if $n\mid \binom{n}{k}$.
%Recently, more attention has been paid to more structured cycles such as tight cycles and loose cycles, and intermediary notions.
Bailey and Stevens~\cite{BS:10} conjectured that $K^k_n$ has a decomposition into tight Hamilton cycles if and only if $n\mid \binom{n}{k}$. This conjecture is generalised in~\cite{KO:14b} to include other notions of cycles such as \defn{loose cycles}.  A related conjecture on wreath decompositions was independently brought forward by Baranyai~\cite{baranyai:79} and Katona in the 1970s.
(If $k$ and $n$ are coprime, then a tight Hamilton cycle coincides with the notion of a \defn{wreath}.)
Approximate results in the sense of packing many edge-disjoint Hamilton cycles into $K^k_n$ have been obtained in~\cite{BF:12,FK:12,FKL:12}.

Some results on Euler tours in hypergraphs have been obtained using the Berge notion (such Euler tours are defined analogously, except that vertices may be repeated).
In~\cite{LN:10}, it is shown that the problem of deciding whether a $k$-graph has a Berge Euler tour is NP-complete for all~$k>2$. On the other hand, a characterization is obtained for so-called `strongly connected' $k$-graphs: such a $k$-graph $G$ has a Berge Euler tour if and only if the number of odd degree vertices of $G$ is at most $(k-2)|E(G)|$. 
%Whilst this condition reduces to Euler's theorem for $k=2$, there is no generic divisibility condition for $k\ge 3$ in the sense that, for example, all $k$-graphs $G$ with $|E(G)|\ge |V(G)|$ trivially satisfy the condition. This is due to the flexibility of the Berge notion. Whilst this can be seen as an advantage (in the sense that more $k$-graphs admit such an Euler tour), it also suggests that the condition that all vertex degrees are divisible by~$k$ more closely resembles the graph property of being Eulerian.
The existence of Berge Euler tours has also been investigated with the host hypergraphs being designs~\cite{HH:14,HH:14b,SW:17}.

It is also natural to seek the above
structures within $k$-graphs of large minimum degree. 
To formalize this, for a set $S\In V(G)$ with $0\le |S|\le k$, we let $d_G(S)$ denote the \defn{degree of $S$ in $G$}, that is, the number of edges which contain~$S$.
We let $\delta(G)$ and $\Delta(G)$ denote the minimum and maximum $(k-1)$-degree of a $k$-graph $G$, respectively, that is, the minimum/maximum value of $d_G(S)$ over all $S\In V(G)$ of size~$k-1$.
R\"odl, Ruci\'nski and Szemer\'edi~\cite{RRS:08} showed that a $k$-graph $G$ on $n$ vertices with $\delta(G) \ge (1/2+o(1))n$ contains a tight Hamilton cycle.
Many related results have been obtained, see e.g.~\cite{zhao:16} for a recent survey.
We pose the following question, which would
show that the degree threshold for a tight Euler tour 
and that for a tight Hamilton cycle coincide asymptotically.

\begin{conj}\label{conj:Dirac}
For all $k>2$ and $\eps>0$, there exists $n_0\in \bN$ such that every $k$-graph $G$ on $n\ge n_0$ vertices with $\delta(G)\ge (1/2+\eps)n$ has a tight Euler tour if all vertex degrees are divisible by~$k$.
\end{conj}

It follows from our Theorem~\ref{thm:typical} that this holds with $1/2+\eps$ being replaced by $1-\eps$ for some small~$\eps$.
The following adaptation of a well-known construction shows that the conjecture would be asymptotically best possible for infinitely many~$n$:
Consider the $k$-graph $G$ with vertex set $V(G)=A\cupdot B$, where $|A|=|B|$ is divisible by~$k$, and all possible edges except those which intersect $B$ in precisely one vertex. 
By removing up to $k-1$ perfect matchings from $G[A]$ and from~$G[B]$, we can ensure that all vertex degrees of the resulting $k$-graph $G'$ are divisible by~$k$. Moreover, $G'$ does not have a tight Euler tour.

\section{Proof} \label{sec:proof}

For a $k$-graph $G$, we let $|G|$ and $e(G)$ denote the number of vertices and edges of $G$, respectively.
Given a (closed) walk $\cW$ in $G$, we will often view $\cW$ as the subgraph $(V(G),E(\cW))$ of~$G$ and accordingly use terminology such as $e(\cW)$ and $\Delta(\cW)$.
 
\subsection{Spanning walk} We call a walk $\cW$ in a $k$-graph $G$ \defn{spanning} if every ordered $(k-1)$-set of vertices appears consecutively in~$\cW$ at least once. One important ingredient of our approach for the proof of Theorem~\ref{thm:typical} is to find a sparse spanning walk in a given $k$-graph~$G$. This spanning walk will form a `backbone' structure to which we will subsequently add smaller closed walks until every edge of $G$ is used exactly once.

We show that such a spanning walk can be obtained randomly, by following a self-avoiding random walk for a suitable number of steps. More precisely, let $G$ be a $k$-graph. We define a simple random process $X=(X_t)_{t\in \bN}$ as follows: Arbitrarily choose distinct starting vertices $x_1,\ldots,x_{k-1}\in V(G)$ and set $X_t:=x_t$ for all $t\in [k-1]$.
Moreover, let $G_{k-1}:=G$.
For all $t\geq k$, proceed as follows.
Among all edges in $G_{t-1}$ that contain the $(k-1)$-set $\{X_{t-k+1},\ldots,X_{t-1}\}$
choose one edge $e$ uniformly at random and let $X_t$ be the vertex in $e\sm \{X_{t-k+1},\ldots,X_{t-1}\}$ and set $G_t:=G_{t-1}- e$. If no such edge is available, then terminate the process and set $X_{t'}:=\emptyset$ and $G_{t'}:=G_{t-1}$ for all $t'\geq t$.

Clearly, this yields a walk in $G$ as long as the process does not terminate.
We write $\cW_t$ for the walk $X_1X_2\dots X_{t'}$ in~$G$, where $t'\le t$ is maximal such that $X_{t'}\neq\emptyset$. Note that $E(\cW_t)=E(G)\sm E(G_t)$.

We will only perform a very crude analysis of this process here, which is sufficient for our purposes. Clearly, for every ordered $(k-1)$-set of the vertices of $G$ to appear in the walk, we need a walk of length $\Omega(|G|^{k-1})$. 
We show that if we follow the random walk for a slightly larger number of steps, then with high probability, every ordered $(k-1)$-set of vertices will indeed appear at least once, and the walk will still be sparse in the sense that no $(k-1)$-set is contained in too many edges of the walk.

For the analysis of the process, we will use the following Chernoff bound. 
It follows directly from the usual Chernoff bound by observing that the moderately dependent Bernoulli variables in our first case
are stochastically dominated by a Binomial random variable with parameters $n$ and $p^+$ (and similarly for the second case).

\begin{lemma}\label{lem:chernoff}
Suppose $X_1,\dots,X_n$ are Bernoulli random variables, and let $X:=\sum_{i=1}^n X_i$. Suppose $0\le \eps \le 3/2$. If for all $i\in[n]$, we have $\prob{X_i=1\mid X_1,\dots,X_{i-1}} \le p^+$, then $$\prob{X \ge (1+\eps) np^+ } \leq \eul^{-\eps^2 np^+/3}.$$ Similarly, if for all $i\in[n]$, we have $\prob{X_i=1\mid X_1,\dots,X_{i-1}} \ge p^-$, then $$\prob{X \le (1-\eps) np^- } \leq \eul^{-\eps^2 np^-/3}.$$
\end{lemma}

The following definition turns out to be a suitable assumption on $G$ which enables a convenient analysis of the process.
We call a $k$-graph $G$ \emph{$\alpha$-connected}
if for all distinct $v_1,\ldots,v_{k-1},v_{k+1},\ldots,v_{2k-1}\in V(G)$,
there exist at least $\alpha |G|$ vertices $v_k$ such that $v_iv_{i+1}\dots v_{i+k-1}\in E(G)$ for all $i\in[k]$. This property is present in natural classes of $k$-graphs. For instance, if $G$ is $(c,k,p)$-typical, then $G$ is $(1-c) p^k$-connected. Similarly, if $\delta(G)\ge (1-\frac{1}{k}+\alpha)|G|$, then $G$ is $k\alpha$-connected.

\begin{lemma} \label{lem:spanning walk}
Let $k\geq 2$ and $\alpha>0$.
Suppose $n$ is sufficiently large in terms of $k$ and~$\alpha$.
Suppose that $G$ is an $\alpha$-connected $k$-graph on $n$ vertices and $X=(X_t)_{t\in \bN}$ is the process defined above. Let $T:=\lflr n^{k-1}\log^2 n\rflr $.
Then with probability at least $1-1/n$, $\cW_T$ is a spanning walk in $G$ and $\Delta(\cW_T)\leq \log^3 n$.
\end{lemma}

\proof
We denote by $\cE_t$ the event that $\Delta(\cW_t)\leq \sqrt{n}$.
We say that distinct $v_1,\ldots,v_{k-1}\in V(G)$ (or simply an ordered $(k-1)$-set) are \emph{met} by $X$ at step~$t\geq k-1$
if $X_{t-k+1+i}=v_{i}$ for all $i\in [k-1]$. A $(k-1)$-set $S$ is \emph{covered} by $X$ at step $t$ if $S\subseteq \Set{X_{t-k+1},\ldots,X_t}$ and $X_t\neq \emptyset$.

In the following, fix distinct $v_1,\ldots,v_{k-1}\in V(G)$.
The key idea is to show that for any given time~$t$, %conditioned on the process $(X_{t'})_{t'\in [t]}$ and the event $\cE_t$,
the probability that $X$ meets or covers $v_1,\ldots, v_{k-1}$ at step $t+2k$ is $\Theta_{\alpha}(n^{-k+1})$.

For $t \geq k-1$,
let $I_t$ be the indicator random variable of the event that $X$ meets $v_1,\ldots,v_{k-1}$ at step $t$,
and let $C_t$ be the indicator random variable of the event that $X$ covers $\Set{v_1,\ldots,v_{k-1}}$ at step~$t$. Note that $d_{\cW_t}(\Set{v_1,\dots,v_{k-1}})= \sum_{t'=k}^t C_{t'}$.

Consider $t\ge k-1$ and suppose we know the outcome of the process up to and including step~$t$. Clearly, in every step of the process, there are at most $n$ choices for the next vertex. Moreover, if $\cE_t$ holds, then $G_{t-1+i}$ will be $2\alpha/3$-connected and $\delta(G_{t-1+i}) \ge 2\alpha n/3$ for all $i\in [2k]$. In particular, in each of the following $2k$ steps, the process has at least $\alpha n/2$ choices for the next vertex $X_{t+i}$ (even if we require that $X_{t+i}\notin \Set{v_1,\ldots,v_{k-1}}$). Thus, the process will not terminate within the next $2k$~steps.

We claim that
\begin{align} \label{eq:lbound}
	\pr[I_{t+2k}=1 \mid (X_{t'})_{t'\in [t]}, \cE_t ] \ge \frac{(\alpha n/2)^{k+1}}{n^{2k}}
	= \frac{\alpha^{k+1}}{2^{k+1}n^{k-1}}=:p^-.
\end{align}
Clearly, there are at most $n^{2k}$ choices for the vertices $X_{t+1},\ldots,X_{t+2k}$.
Moreover, for at least $(\alpha n/2)^{k+1}$ choices of $X_{t+1},\ldots,X_{t+2k}$, the vertices $v_1,\ldots,v_{k-1}$ are met at step $t+2k$. This is
because for $X_{t+1},\ldots,X_{t+k}$ there are at least $(\alpha n/2)^k$ choices that avoid $v_1,\ldots,v_{k-1}$,
and then there are at least $\alpha n/2$ choices for $X_{t+k+1}$ such that $X_{t+k+1+i}=v_i$ for all $i\in [k-1]$ is a valid choice for the process. (Here, the step of choosing $X_{t+k+1}$ is the part where the definition of $\alpha$-connectedness is crucial.)

We also claim that
\begin{align} \label{eq:ubound}
	\pr[C_{t+2k} =1 \mid (X_{t'})_{t'\in [t]}, \cE_t] \le \frac{k!n^{k+1}}{(\alpha n/2)^{2k}}
	= \frac{2^{2k}k!}{\alpha^{2k}n^{k-1}}=:p^+.
\end{align}
To prove this claim, we make three observations.
Firstly, recall that in each of the next $2k$~steps, the process has at least $\alpha n/2$ choices for the next vertex.
Secondly, note that if $C_{t+2k}=1$, then $\Set{X_{t+k+1},\ldots,X_{t+2k}}\supseteq \Set{v_1,\dots,v_{k-1}}$. There are at most $k!$ ways of assigning $v_1,\ldots,v_{k-1}$ to $X_{t+k+1},\ldots,X_{t+2k}$. Thirdly, there are at most $n^{k+1}$ choices for the remaining $k+1$ vertices.

Note that the probability estimates \eqref{eq:lbound} and \eqref{eq:ubound} rely on the assumption that $\cE_t$ holds. To account for the complementary case, we define auxiliary $0/1$ random variables $Y_t^-$ and $Y_t^+$ for $t\ge 4k+ 1$ as follows.
Let $Y_t^- :=I_t$ and $Y_t^+ :=C_t$ if $\cE_{t-2k}$ holds and otherwise let $Y_t^- :=1$ with probability $p^-$
and $Y_t^+ :=1$ with probability $p^+$ independently of all other random choices.

Since the bounds \eqref{eq:lbound} and \eqref{eq:ubound} only hold if we condition on the process until $2k$ steps earlier, we consider $2k$ disjoint subsequences and analyse each subsequence individually.
Define $T':= \lfloor T/2k \rfloor$ and for all $i\in [2k]$, define $Z_{i}^\pm :=\sum_{t'=1}^{T'}{Y_{i+2k(t'+1)}^\pm}$. %and $Z_{i}^+ :=\sum_{t= 2}^{T'}{Y_{2kt+i}^+}$.
Observe that for each $i\in [2k]$ and all $t'\in \bN$, we have
\begin{align*}%\label{eq:lbound}
	\pr[Y_{i+2k(t'+1)}^-=1 \mid Y_{i+4k}^-,Y_{i+6k}^-,\dots, Y_{i+2kt'}^-] \ge p^-, \\
	\pr[Y_{i+2k(t'+1)}^+=1 \mid Y_{i+4k}^+,Y_{i+6k}^+,\dots, Y_{i+2kt'}^+] \le p^+.
\end{align*}
This follows from \eqref{eq:lbound} and \eqref{eq:ubound} if $\cE_{i+2kt'}$ holds and from the definition of $Y_{i+2k(t'+1)}^\pm$ otherwise.
%Thus, $Z_{i}^-$ stochastically dominates a binomial random variable with parameters $T'$ and $p^-$.
%Similarly, $Z_{i}^+$ is stochastically dominated by a binomial random variable with parameters $T'$ and $p^+$.
Note that $T'\cdot p^- \geq \log^{3/2} n$ and $\log^{2} n\le T'\cdot p^+ \leq \log^{5/2} n$.
Hence, by Lemma~\ref{lem:chernoff} we conclude that
\begin{align*}
	\pr[Z_{i}^-\leq (\log^{3/2} n)/2] \le \eul^{-(\log^{3/2} n)/12} \text{ and } \pr[Z_{i}^+\geq 2\log^{5/2} n] \le \eul^{-(\log^2 n)/3}.
\end{align*}
Therefore, with probability at least $1-1/n$ say, we have $Z_{i}^-\geq (\log^{3/2} n)/2$ and $Z_{i}^+\leq 2\log^{5/2} n$ for all $i\in [2k]$
and all choices of $v_1,\ldots,v_{k-1}$ simultaneously.

Finally, suppose that the above hold. We claim that $\cE_T$ holds. Suppose not. Let $t_0\le T$ be minimal such that $d_{\cW_{t_0}}(S) > \sqrt{n}$ for some $(k-1)$-set~$S$. Then $Y_t^+=C_t$ for all $t\leq t_0$, and hence $d_{\cW_{t_0}}(S)\le 5k + \sum_{i=1}^{2k}Z_i^+$ (the first term accounts for the roughly $4k$ steps that are not taken into account by the $Z_i^+$ variables), a contradiction.
This implies that $Y_t^-=I_t$ and $Y_t^+=C_t$ for all $t\leq T$.
Consequently, all ordered $(k-1)$-sets are met at least once until step $T$ and $\Delta(\cW_T)\leq \log^3 n$.
\endproof

\subsection{$F$-decompositions}

In 1976, Wilson~\cite{wilson:76} proved the fundamental result that given any graph~$F$, for sufficiently large~$n$, the complete graph $K_n$ has an $F$-decomposition whenever it satisfies some necessary divisibility conditions (see below).
This was generalised to hypergraphs in~\cite{GKLO:17}. In order to formally state the required result, we define the following. Let $G$ and $F$ be $k$-graphs, where $F$ is non-empty. An \defn{$F$-decomposition of $G$} is a collection of copies of $F$ in $G$ such that every edge of $G$ is contained in exactly one of these copies. It is easy to see that the existence of an $F$-decomposition necessitates certain divisibility conditions. For instance, we surely need $e(F)\mid e(G)$. More generally, define $d_F(i):=\gcd\set{d_F(S)}{S\in \binom{V(F)}{i}}$ for all $0\le i\le k-1$.\COMMENT{As long as $F$ is not edgeless, this is well-defined.} Note that $d_F(0)=e(F)$.
%So if $F$ is the Fano plane, we have $Deg(F)=(7,3,1)$.
Now, $G$ is called \defn{$F$-divisible} if $d_F(i)\mid d_G(S)$ for all $0\le i\le k-1$ and all $S\in \binom{V(G)}{i}$.
It is easy to see that $G$ must be $F$-divisible in order to admit an $F$-decomposition. The converse implication is in general not true. However, if $G$ is a large typical $k$-graph, then divisibility guarantees the existence of a decomposition. For $G=K^k_n$, this generalises Wilson's theorem to hypergraphs.

\begin{theorem}[\cite{GKLO:17}] \label{thm:design}
For all $k\in \bN$, $p\in[0,1]$ and any $k$-graph~$F$,
there exist $c>0$ and $h,n_0\in\bN$ such that the following holds.
Suppose that $G$ is a $(c,h,p)$-typical $k$-graph on at least $n_0$~vertices. Then $G$ has an $F$-decomposition whenever it is $F$-divisible.
\end{theorem}
We remark that explicit bounds for $c$ and $h$ were obtained in~\cite{GKLO:17}. Using these one can also obtain such explicit bounds in Theorem~\ref{thm:typical}. For a subsequent alternative proof of Theorem~\ref{thm:design} see~\cite{keevash:18b}.

%\begin{theorem}[\cite{GKLO:17}] \label{thm:design}
%For all $f,k\in \bN$ with $f>k$ and all $c,p\in (0,1]$ with
%\begin{align*}
%c \le 0.9(p/2)^{h}/(q^k4^q), \mbox{ where } q:=2f\cdot f! \mbox{ and } h:=2^k\binom{q+k}{k},
%\end{align*}
%there exists $n_0\in \bN$ such that the following holds for all $n\ge n_0$.
%Let $F$ be any $k$-graph on $f$~vertices. Suppose that $G$ is a $(c,h,p)$-typical $k$-graph on $n$~vertices. Then $G$ has an $F$-decomposition if it is $F$-divisible.
%\end{theorem}

Let $C^{k}_\ell$ denote the tight $k$-uniform cycle of length $\ell$, that is, the vertices of $C^{k}_\ell$ are $v_1,\dots,v_\ell$, and the edges are all the $k$-tuples of the form $\Set{v_i,v_{i+1},\dots,v_{i+k-1}}$, with indices modulo~$\ell$.

Here, we will apply Theorem~\ref{thm:design} with $F=C^{k}_{2k}$. Clearly, we have $d_{C^{k}_{2k}}(0)=e(C^{k}_{2k})=2k$ and $d_{C^{k}_{2k}}(1)=k$. Moreover, for every $i\in \Set{2,\dots,k-1}$, we have $d_{C^{k}_{2k}}(\Set{v_1,\dots,v_{i-1},v_k})=1$ and hence $d_{C^{k}_{2k}}(i)=1$. Conveniently, a $k$-graph $G$ is thus $C^{k}_{2k}$-divisible whenever $2k\mid e(G)$ and $k \mid d_G(v)$ for all $v\in V(G)$.

\subsection{Proof of Theorem~\ref{thm:typical}}
We can now prove our main theorem.

\lateproof{Theorem~\ref{thm:typical}}
Given $k$ and $p$, choose $c>0$ sufficiently small and $h,n_0$ sufficiently large. In particular, we assume that $h\ge k$ and that we can apply Theorem~\ref{thm:design} with $k$, $p$, $2c$, $h$, $n_0$, $C^{k}_{2k}$ playing the roles of $k,p,c,h,n_0,F$.
Suppose that $G$ is a $(c,h,p)$-typical $k$-graph on $n\ge n_0$ vertices with all vertex degrees divisible by~$k$. Since $G$ is $(1-c)p^k$-connected, by Lemma~\ref{lem:spanning walk}, there exists a spanning walk $\cW=v_1v_2\dots v_\ell$ in~$G$ such that $\Delta(\cW)\le \log^3 n$. Next, we extend $\cW$ to a closed walk~$\cW'$. Choose $k\le \ell'\le 3k-1$ such that $\ell'\equiv e(G)-\ell \mod{2k}$. Now, find distinct vertices $v_{\ell+1},\dots,v_{\ell+\ell'}\in V(G)\sm \Set{v_1,\dots,v_{k-1},v_{\ell-k+2},\dots,v_\ell}$ such that $v_iv_{i+1}\dots v_{i+k-1}\in E(G)\sm E(\cW)$ for all $i$ with $\ell-k+1<i\le \ell+\ell'$, with indices modulo~$\ell+\ell'$. We can find such vertices one-by-one using the typicality of~$G$. Indeed, when finding $v_{\ell+j}$, we need to ensure that $v_{\ell+j}$ belongs to the neighbourhood of (at most) $k$ specific $(k-1)$-subsets $S\In V(G)$, but does not form an edge of $\cW$ with any of these. Since $\Delta(\cW)\le \log^3 n$, there are at least $(1-c)p^k n - k\log^3 n\ge p^k n/2$ such vertices, from which we can choose one whilst also avoiding previously chosen vertices and $\Set{v_1,\dots,v_{k-1},v_{\ell-k+2},\dots,v_\ell}$.
Let $\cW':=v_1\dots v_{\ell+\ell'}$. Clearly, $\cW'$ is a spanning closed walk in~$G$ and $\Delta(\cW')\le 2\log^3 n$.

Now, let $G':=G-E(\cW')$. We have $e(G') =e(G)-(\ell+\ell') \equiv 0\mod{2k}$. Moreover, since $\cW'$ is a closed walk, we have $k\mid d_{\cW'}(v)$ for all $v\in V(G)$. Combining this with the initial divisibility condition of $G$, we have that $k\mid d_{G'}(v)$ for all $v\in V(G')$. Thus, $G'$ is $C^{k}_{2k}$-divisible.
Moreover, since $\Delta(\cW')\le 2\log^3 n $, we have that $G'$ is $(2c,h,p)$-typical.\COMMENT{If $H$ is $(c,h,p)$-typical and $\Delta(L)\le \gamma n$ with $V(L)=V(H)$, then $H\bigtriangleup L$ is $(c+hp^{-h}\gamma, h,p)$-typical.}  Invoking Theorem~\ref{thm:design}, we conclude that $G'$ has a $C^{k}_{2k}$-decomposition~$\cC$. We can now simply incorporate each cycle of $\cC$ one-by-one into the spanning closed walk~$\cW'$. For this, suppose that $\cW''$ is the current spanning closed walk and let $C\in \cC$ be a copy of $C^{k}_{2k}$ with vertices $v_1,\dots,v_{2k}$ appearing in this order on~$C$. Since $\cW''$ is spanning, $v_1,\dots,v_{k-1}$ appear consecutively in $\cW''$, say $\cW''=\cW''_1 v_1\dots v_{k-1} \cW_2''$. We can then simply replace $\cW''$ with $\cW''_1 v_1\dots v_{2k} v_1\dots v_{k-1} \cW_2''$ to obtain a new spanning closed walk $\cW'''$ with $E(\cW''')=E(\cW'')\cup E(C)$.
Adding all cycles of $\cC$ in this way yields the desired Euler tour.
\endproof

As pointed out after Theorem~\ref{thm:typical}, our proof can also be easily adapted to obtain a tight Euler trail. If a ($2$-)graph has an Euler trail, then there are precisely two vertices of odd degree. If an Euler trail in a $k$-graph $G$ starts with the sequence $v_1\dots v_{k-1}$ and ends with the sequence $w_{k-1}\dots w_1$, then (assuming that $v_1,\dots,v_{k-1},w_1,\dots,w_{k-1}$ are distinct) we must have $d_G(v_i),d_G(w_i)\equiv i\mod{k}$ for all $i\in[k-1]$, and all other vertex degrees are divisible by~$k$.
On the other hand, if these conditions hold and $G$ is typical, then $G$ has an Euler trail. For this, in the above proof of Theorem~\ref{thm:typical}, instead of extending the spanning walk $\cW$ to a closed walk, one simply extends both ends to the designated start and end $(k-1)$-tuples.

\providecommand{\bysame}{\leavevmode\hbox to3em{\hrulefill}\thinspace}
\providecommand{\MR}{\relax\ifhmode\unskip\space\fi MR }
\providecommand{\MRhref}[2]{%
  \href{http://www.ams.org/mathscinet-getitem?mr=#1}{#2}
}
\providecommand{\href}[2]{#2}

%\bibliographystyle{../../amsplain_v2.0customized}
%\bibliography{../../References}

\end{document}